\documentclass[12pt]{article}

\usepackage{amssymb}
\usepackage{amsmath}
\usepackage{amsfonts}
\usepackage{graphicx}

\def\bea{\begin{eqnarray}}
\def\eea{\end{eqnarray}}
\newcommand{\la}{\label}
\newcommand{\nn}{\nonumber}

\font\mybb=msbm10 at 11pt

\def\bb#1{\hbox{\mybb#1}}

\def\bZ {\bb{Z}}
\def\bR {\bb{R}}

\def\bC {\bb{C}}

\def\lc{\lrcorner}
\def\cl{{\mathrm {Cliff}}}

\newcommand{\tM}{\text{\tiny $M$}}
\newcommand{\tB}{\text{\tiny $B$}}
\newcommand{\tA}{\text{\tiny $A$}}

\newcommand{\tN}{\text{\tiny $N$}}
\newcommand{\tL}{\text{\tiny $L$}}
\newcommand{\tP}{\text{\tiny $P$}}
\newcommand{\tR}{\text{\tiny $R$}}
\newcommand{\tQ}{\text{\tiny $Q$}}

\def\e  {\epsilon}
\def\a {\alpha}

\def\sp {{\cal M}}

\begin{document}

\begin{titlepage}
\begin{center}
\vspace*{3.50cm}
{\Large \bf Solution of heterotic Killing spinor equations and special geometry}\\[.2cm]

\vspace{1.5cm}
 {\large  George Papadopoulos$^1$ and Ulf Gran$^2$}

\vspace{1.5cm}

${}^1$ Department of Mathematics\\
King's College London\\
Strand\\
London WC2R 2LS, UK\\

\vspace{0.5cm}
${}^2$ Fundamental Physics\\
Chalmers University of Technology\\
SE-412 96 G\"oteborg, Sweden\\

\end{center}

\vskip 1.5 cm
\begin{abstract}
We outline  the solution of the Killing spinor equations of the heterotic supergravity. In addition,
we describe the classification of all half supersymmetric solutions.

\end{abstract}
\end{titlepage}

\section{Introduction}

Supersymmetric supergravity backgrounds are solutions of the field equations of supergravity theories which in addition
solve a set of first order equations, the Killing spinor equations. These solutions are triplets
$({\sp},g,F)$, where $\sp$ is a Lorentzian manifold with metric $g$, and $F$ are the fluxes of supergravity theories which
is a collection of forms on $\sp$.
The field equations of supergravity theories consist of the Einstein equation as well as appropriate Maxwell type of  equations
for $F$. The Killing spinor equations   are determined  from the supersymmetry transformations of the fermions\footnote{After considering
the supersymmetry transformations, in what follows all the fermionic fields are set to zero.} of the supergravity
theories. Moreover their integrability conditions imply some of the supergravity field equations.

Recently, there is much interest in systematically understanding the supersymmetric solutions of the supergravity theories.
This has been mostly motivated by the applications that these solutions have in string theory, M-theory and in the AdS/CFT correspondence. Apart from
this, the supersymmetric supergravity solutions are the gravitational analogues of gauge theory solitons and instantons, and so
their classification is interesting in its own right.

The main aim of this article is to outline the classification of the solutions of the Killing spinor equations
of the heterotic supergravity \cite{het1, het2, het3}. Moreover, all supersymmetric solutions which preserve 8 Killing spinors will be described \cite{het4}.
This material is partly based  on work done in collaboration with Diederik Roest, Philipp Lohrmann and Peter Sloane
as well as on material published by one of the authors.
In addition,  this paper  contains a refinement of the results
of the first two papers. In particular, a more concise description of the geometry of
the backgrounds with non-compact holonomy is given in terms of
certain  Clifford algebras of endomorphisms.

This paper is organized as follows. In section two, the Killing spinor and field equations of the heterotic
supergravity are given. We also summarize the main ingredients of the method  that we use to solve the Killing spinor equations.
In section three, the gravitino Killing spinor equation is solved. In section four, an outline of the solution
of the dilatino Killing spinor equation is given. In section five,
the geometry of supersymmetric backgrounds with non-compact and compact holonomy is described. In section six,
we solve the field equations of the heterotic supergravity   for all backgrounds preserving 8 supersymmetries, i.e.~we describe all half supersymmetric backgrounds.

\section{Killing spinor and field equations}

\subsection{Killing spinor and field equations}

The spacetime is a 10-dimensional Lorentzian manifold $\sp$. The bosonic fields of heterotic supergravity   are a metric $g$,
a 3-form field strength $H$, the dilaton scalar
$\Phi$,  and a gauge connection $A$ with curvature $F=dA- A\wedge A$. The gauge group of $A$ is either $E_8\times E_8$ or
$Spin(32)/\bZ_2$. Though this restriction on the gauge group does not affect most of the analysis that will follow.

 The gravitino, gaugino and dilatino Killing spinor
equations of the heterotic supergravity are
\bea
 {\cal D}_\tM \e &\equiv& \hat\nabla_\tM\e+{\cal O}(\a'^2)=0~,~~~
 {\cal A} \e \equiv (\Gamma^\tM\partial_\tM\Phi-{1\over12} H_{\tM\tN\tL}\Gamma^{\tM\tN\tL})\e+{\cal O}(\a'^2)=0~,~~~
 \cr
  {\cal F} \e  &\equiv& F_{\tM\tN}\Gamma^{\tM\tN} \e+{\cal O}(\a'^2)=0~,~~~
\la{kse}
\eea
respectively, where $\epsilon$ is a real positive chirality  spinor (Majorana-Weyl) of $Spin(9,1)$ and $\hat\nabla=\nabla+{1\over2} H$
is a metric connection with torsion $H$. Moreover,  $\{\Gamma_\tM\}$ is a basis of the Clifford algebra ${\rm Cliff}(\bR^{9,1})$,
\bea
\Gamma_{\tM} \Gamma_{\tN}+
\Gamma_{\tN} \Gamma_{\tM}=2 g_{\tM\tN}~,
\eea
 and $M,N,L=0,\dots,9$.
More details about the notation can be found in \cite{het1, het2}. The Killing spinor equations have been expressed as an expansion in the parameter $\a'$. They are known to the order indicated
but it is expected that they receive  corrections to higher orders.

The 3-form field strength $H$ is not closed but is modified at order $\a'$ because of the Green-Schwarz anomaly cancelation mechanism as
\bea
dH=- { \a'\over 4}\big( {\rm tr} \,\,\check R^2- {\rm tr}\,\, F^2\big)+{\cal O}(\a'^2)~,
\la{anb}
\eea
where $\check R$ is the curvature of $\check\nabla=\nabla-{1\over2} H$.

 The field equations (in the string frame) to lowest
order in $\alpha'$  are
\bea
E_{\tM\tN}&\equiv&R_{\tM\tN}+{1\over4} H^\tR{}_{\tM\tL} H^\tL{}_{\tN\tR}+2\nabla_\tM\partial_\tN\Phi
\cr&&~~~~~~~~~~+{\a'\over4} [\check R_{\tM\tL,\tQ\tR} \check R_\tN{}^{\tL,\tQ\tR}-F_{\tM\tL ab} F_\tN{}^{\tL ab}]+{\cal O}(\a'^2)=0~,
\cr
LH_{\tP\tR}&\equiv&\nabla_\tM[e^{-2\Phi} H^\tM{}_{\tP\tR}]+{\cal O}(\a'^2)=0~,
\cr
LF_\tM&\equiv&\hat\nabla^\tM[e^{-2\Phi} F_{\tM\tN}]+{\cal O}(\a'^2)=0~.
\la{feqn}
\eea
The linear term in $\a'$ in the Einstein equation, which arises from the 2-loop sigma model beta function calculation \cite{town}, is
necessary for consistency with (\ref{anb}), see e.g.~\cite{tsimpis}. The remaining two field equations are Maxwell type of equations
for the 3-form flux $H$ and the 2-form gauge field strength $F$.
The field equation for the dilaton is implied from those above up to a constant.

\subsection{Method}

The method we shall use to solve the Killing spinor equations of the heterotic supergravity is {\it spinorial geometry} \cite{uggp}.
It is based on
\begin{itemize}

\item the gauge symmetry of Killing spinor equations,

\item a description of spinors in term of forms,

\item a harmonic oscillator basis in the space of spinors.

\end{itemize}
The basic strategy is to use the gauge symmetry of the Killing spinor equations  to choose a canonical form for the Killing spinors or their normals.
Then writing the Killing spinors  in terms of forms, these can be substituted into the Killing spinor equations.
The resulting expressions are solved by utilizing the linearity of the Killing spinor equations and expanding them in the harmonic oscillator
 basis in the space of spinors.

The above method is very effective particularly for the solutions of the Killing spinor equations with small
or near maximal number of supersymmetries\footnote{For the classification of  near maximally supersymmetric
backgrounds  using  spinorial geometry see \cite{n31}.}. It can be  implemented equally efficiently in analytic or
computer calculations.

Returning to the Killing spinor equations of heterotic supergravity, it is convenient to solve them in the order
\bea
{\rm gravitino}\rightarrow {\rm gaugino}\rightarrow {\rm dilatino}~.
\nn
\eea
The solution of the gaugino Killing spinor equation has been given in \cite{het3}, and it is similar
to that of the gravitino. Because of this in the analysis that follows, we shall focus on the
solution of the gravitino and dilatino Killing spinor equations \cite{het1, het2}.

To apply the spinorial geometry method  to the heterotic supergravity, first observe that the gauge symmetry of the
Killing spinor equations (\ref{kse}) is $Spin(9,1)$. This coincides with the holonomy group of $\hat\nabla$, ${\rm hol}(\hat\nabla)$,
for  generic backgrounds. This equality is the main reason that all the solutions
of the Killing spinor equations of the heterotic supergravity can be found.

\subsection{Spinors}\la{spinors}

One of the ingredients of  spinorial geometry is the description of spinors in terms of forms. This
is a well-known realization of the spinor representations, see  e.g.~\cite{ harvey}, and it has been used in \cite{wang} to give
explicitly the parallel spinors of Riemannian manifolds with special holonomy.
This description of spinors can be extended to the Lorentzian case. For later use,  we  give  the form realization of
 spinor representations\footnote{The spin groups considered here are the double covers of the component of the Lorentz group connected to the identity.} of $Spin(9,1)$, see also \cite{het1}.

Consider  $\bC^5=\bC<e^1,\dots,e^5>$, where  $e^1,\dots,e^5$ is a Hermitian basis with respect to the $<\cdot, \cdot>$ inner product.
The space of Dirac spinors of $Spin(9,1)$ is
$\Delta_c=\Lambda^*(\bC^5)$.
The basis $\{\Gamma_\tA\}$ of Clifford algebra ${\rm Cliff}(\bR^{9,1})$  acts on  $\Delta_c$ as
\bea
\Gamma_0\psi&=& -e_5\wedge\psi +e_5\lc\psi~,~~~~ \Gamma_5\psi=
e_5\wedge\psi+e_5\lc \psi~,
\cr
\Gamma_i\psi&=& e_i\wedge \psi+ e_i\lc \psi~,~~~~~~
\Gamma_{5+i}\psi= i e_i\wedge\psi-ie_i\lc\psi~, ~~~~~~i=1,\dots,4,
\eea
where $\psi\in \Delta_c$ and $\lc$ is the adjoint operation of $\wedge$ with respect to $<\cdot,  \cdot>$.
It is easy to verify that $\{\Gamma_\tA\}$ satisfies the Clifford algebra relation $\Gamma_\tA \Gamma_\tB+
\Gamma_\tB \Gamma_\tA=2 \eta_{\tA\tB}$, where $\eta$ is the Minkowski metric.
$\Delta_c$ is a reducible $Spin(9,1)$ representation and decomposes into two complex chiral representations
 $\Delta_c^+=\Lambda^{{\rm even}}( \bC^5)$
and  $\Delta_c^-=\Lambda^{{\rm odd}}(\bC^5)$. These are the complex Weyl representations
of $Spin(9,1)$.

It is well-known that $Spin(9,1)$ admits two inequivalent {\it real} chiral representations, the Majorana-Weyl representations.
These are constructed by imposing a reality condition on $\Delta^\pm_c$. This is achieved by using the reality map $R=\Gamma_{6789}*$
which is anti-linear, $R^2=1$, and commutes in the action of $Spin(9,1)$. So the real spinors satisfy
\bea
\eta^*=\Gamma_{6789}\eta~.
\eea
  For example the real and imaginary components of the complex spinor $1$ are
$1+e_{1234}$ and $i ( 1- e_{1234})$, respectively, where
$e_{1234}=e_1\wedge e_2\wedge e_3\wedge e_4$. We denote the real subspaces
of $\Delta_c^\pm$ with $\Delta^\pm_{16}$.

The spacetime form bilinears associated with the spinors $\psi,\theta$.
 are given as
\bea
\alpha(\psi, \theta)\equiv{1\over k!} B(\psi,\Gamma_{A_1\dots A_k} \theta)\,\,
e^{A_1}\wedge\dots\wedge e^{A_k}~,~~~~~~~k=0,\dots, 9~.
\la{forms}
\eea
where
\bea
B(\psi,\theta)= <{\rm B}(\psi^*), \theta>~,~~~~~~~~
\eea
 is
the $Spin(9,1)$-invariant Majorana bilinear inner product on $\Delta_c$  and
the linear map ${\rm B}$ is ${\rm B}=\Gamma_{06789}$.

\section{Gravitino Killing spinor equation} \label{sc.b}
Let us assume that the spacetime $\sp$ is simply connected. To investigate the solutions of the gravitino Killing spinor
equation \cite{het1, het2}, consider the integrability condition
\bea
\hat R\,\e=0~.
\eea
This equation has  solutions if either $\hat R=0$, and so $\sp$ is parallelisable, or the solutions  $(\e_1, \dots, \e_\tL)$ have a non-trivial
isotropy group ${\rm Stab}(\e_1, \dots, \e_\tL)\subset Spin(9,1)$ and
\bea
{\rm hol}(\hat\nabla)\subseteq {\rm Stab}(\e_1, \dots, \e_\tL)~.
\eea
In the former case, $\sp$ is either a Lorentzian manifold or a product of a Lorentzian group manifold with $S^7$ \cite{jose1}. If in addition,
one assumes that $dH=0$, then $\sp$ is a Lorentzian group manifold. The Lorentzian groups manifolds have been classified in \cite{medina}.
Locally up to dimension ten, they are products of the Lorentzian groups $\bR$, $SL(2,\bR)$, $CW_{2k}$ with the Riemannian  group manifolds
$U(1)$, $SU(2)$ and $SU(3)$, where $CW_{2k}$ are the group manifolds\footnote{These are
plane waves with wave profile given by the square of a skew-symmetric matrix.} associated with the Cahen-Wallace spaces.

In the latter case, one can determine the subgroups of $Spin(9,1)$ which are isotropy groups of spinors. These  have been tabulated in table 1.
This table has been
constructed in stages  \cite{jose2,  het1, het2}.
In the same table, a basis in the space of parallel spinors is given for each case. These bases have been written down
explicitly using the form notation for spinors explained  in  section \ref{spinors} and they are determined up to  $Spin(9,1)$ gauge transformations.

A straightforward observation of results tabulated in table 1 reveals that there are two types of isotropy groups of
spinors that occur distinguished by their topology,
the compact and non-compact ones. The non-compact isotropy groups
are of the type $K\times\bR^8$, where $K$ is compact. As we shall explain this distinction is useful in the description of
geometry of the associated spacetimes.  Most of the isotropy groups that occur are of the Berger type.
However there are some exceptions which do not appear in the Berger list. These are whenever ${\rm Stab}(\e_1, \dots, \e_\tL)$ is
$\times^2 SU(2)\ltimes \bR^8$, $SU(2)\ltimes \bR^8$,
$U(1)\ltimes \bR^8$ and $\bR^8$.

\begin{table}[ht]
\begin{center}
\begin{tabular}{|c|c|c|c|}\hline
$L$ & ${\mathrm Stab}(\e_1,\dots,\e_L)$ & $\Sigma({\cal P})$& $\e_1,\dots,\e_L$
 \\ \hline \hline
$1$ & $Spin(7)\ltimes\bR^8$& $Spin(1,1)$ &$1+e_{1234}$\\
\hline
$2$ &  $SU(4)\ltimes\bR^8$&$Spin(1,1)\times U(1)$& $1$
\\ \hline
$3$ & $Sp(2)\ltimes\bR^8$&$Spin(1,1)\times SU(2)$ &$1,~~ i(e_{12}+e_{34})$
\\ \hline
$4$ & $\times^2 SU(2)\ltimes\bR^8$&$Spin(1,1)\times(\times^2 Sp(1))$& $1,~~e_{12}$
\\ \hline
$5$ & $SU(2)\ltimes\bR^8$&$Spin(1,1)\times Sp(2)$&$1,~~ e_{12}, ~~e_{13}+e_{24}$
\\ \hline
$6$ & $U(1)\ltimes\bR^8$&$Spin(1,1)\times SU(4)$& $1, ~~e_{12},~~ e_{13}$
\\ \hline
$8$ & $\bR^8$&$Spin(1,1)\times Spin(8)$& $1, e_{ij}, i,j\leq 4$
\\ \hline \hline
$2$ & $G_2$&$Spin(2,1)$&$1+e_{1234}, e_{15}+e_{2345}$
\\ \hline
$4$ & $SU(3)$&$Spin(3,1)\times U(1)$& $1, e_{15}$
\\ \hline
$8$ & $SU(2)$&$Spin(5,1)\times SU(2)$& $1, e_{12}, e_{15}, e_{25}$
\\ \hline
$16$ & $\{1\}$&$Spin(9,1)$&$ 1, e_{ij}, e_{i5}$
\\ \hline
\end{tabular}
\end{center}
\caption{In the  columns are the numbers of parallel spinors, their isotropy groups  and their $\Sigma({\cal P})$ groups, respectively.
The $\Sigma({\cal P})$
groups are a product of a $Spin$ group and an R-symmetry group of a lower-dimensional supergravity theory.}
\end{table}

\section{Dilatino Killing spinor equation} \label{}

Suppose that we have a solution of the gravitino Killing spinor equation and the $\hat\nabla$-parallel spinors span an L-plane
${\cal P}_\tL$. Typically only some of the $\hat\nabla$-parallel spinors will be Killing, i.e.~they will  solve both the gravitino and
 dilatino Killing spinor
equations. Following \cite{het2} to solve the dilatino Killing spinor equation, one has to choose representatives for the Killing spinors
up to $Spin(9,1)$ gauge transformations. It turns out that given the $\hat\nabla$-parallel spinors, a suitable choice of gauge
transformations is
\bea
\Sigma({\cal P}_\tL)={\rm Stab}({\cal P}_\tL)/{\rm Stab}(\e_1,\dots, \e_\tL)~,
\eea
where ${\rm Stab}({\cal P}_\tL)=\{\ell\in Spin(9,1)\vert~~ \ell{\cal P}_\tL\subseteq {\cal P}_\tL\}$. The quotient with ${\rm Stab}(\e_1,\dots, \e_\tL)$
is taken because this subgroup acts with the identity on ${\cal P}_\tL$.
The $\Sigma({\cal P}_\tL)$ groups have been tabulated in table 1.

The analysis of the solutions of the dilatino Killing spinor equation proceeds as follows. Given a solution of the
gravitino Killing spinor equations, one determines ${\cal P}_\tL$. Now suppose  only one of the parallel spinors is Killing. This
can be chosen up to $\Sigma({\cal P}_\tL)$ gauge transformations. Therefore the distinct solutions of the dilatino Killing spinor
equation are labeled by the different type of orbits, ${\mathcal O}_{\Sigma({\cal P}_\tL)}({\cal P}_\tL)$
of  $\Sigma({\cal P}_\tL)$ in ${\cal P}_\tL$.

Having established a procedure to choose the first Killing spinor, one can proceed inductively. Let ${\cal K}_\tN$
denote the subspace of ${\cal P}_\tL$ spanned by the first $N$ Killing spinors, $N<L$. One writes
\bea
0\rightarrow {\cal K}_\tN\rightarrow {\cal P}_\tL\rightarrow {\cal P}_\tL/{\cal K}_\tN\rightarrow 0~.
\eea
The task is to determine $ {\cal K}_{\tN+1}$. For this one has to choose an additional Killing spinor $\e_{\tN+1}\in {\cal P}_\tL$ which
is linearly independent from those in ${\cal K}_\tN$. For this, one uses as a gauge group
\bea
{\rm Stab}({\cal K}_\tN)=\{\ell\in \Sigma({\cal P}_\tL)\vert~~~\ell {\cal K}_\tN\subseteq {\cal K}_\tN\}~.
\eea
Because of the linearity of the dilatino Killing spinor equations, one can view the additional Killing  spinor $\e_{\tN+1}$ as element of
${\cal P}_\tL/{\cal K}_\tN$. Thus, the distinct choices of $\e_{\tN+1}$ are labeled by the different type of orbits,
${\mathcal O}_{{\rm Stab}({\cal K}_\tN)}({\cal P}_\tL/{\cal K}_\tN )$, of ${\rm Stab}({\cal K}_\tN)$ in ${\cal P}_\tL/{\cal K}_\tN$.
In practise this procedure is carried out for $N\leq L/2$. If $N>L/2$, then a similar procedure can be devised for selecting
the normals to the Killing spinors in ${\cal P}_\tL$.
Using the above procedure, the dilatino Killing spinor equation has been solved in all cases and the possibilities that arise have been
tabulated in table 2.

\begin{table}
\begin{center}
\begin{tabular}{|c|c|c|}\hline
$L$&   ${\rm Stab}(\e_1, \dots, \e_L)$ &$N$
 \\\hline\hline
 $1$& $Spin(7)\ltimes\bR^8$& ${ 1(1)}$ \\
 \hline
$2$&$SU(4)\ltimes\bR^8$&$1(1),\,\, { 2(1)}$
\\
\hline
$3$&$Sp(2)\ltimes\bR^8$&$1(1),\,\, 2(1),\,\,  { 3(1)}$
\\
\hline
$4$&$(\times^2SU(2))\ltimes\bR^8$&$1(1),\,\, 2(1),\,\, 3(1), \,\, { 4(1)}$
\\
\hline
$5$&$SU(2)\ltimes\bR^8$&$1(1),\,\, 2(1),\,\, 3(1),\,\, 4(1),\,\, { 5(1)}$
\\
\hline
$6$&$U(1)\ltimes\bR^8$&$1(1),\,\, 2(1),\,\, 3(1),\,\, 4(1),\,\, 5(1),\,\,  { 6(1)}$
\\
\hline
$8$&$\bR^8$&$1(1),\,\, 2(1),\,\, 3(1),\,\, 4(1),\,\, 5(1),\,\, 6(1),\,\, { 7(1)},\,\,  { 8(1)}$
\\
\hline \hline
$2$&$G_2$& $1(1),\,\,{ 2(1)}$
\\
\hline
$4$&$SU(3)$&$1(1),\,\, 2(2), \,\,3(1),\,\, { 4(1)}$
\\
\hline
$8$&$SU(2)$&$1(1),\,\, 2(2),\,\, 3(3),\,\, 4(6),\,\, 5(3),\,\, 6(2),\,\, { 7(1)},\,\,{ 8(1)}$
\\
\hline
$16$&$\{1\}$&$8(2),\,\, 10(1),\,\, 12(1),\,\, 14(1), \,\,{ 16(1)}$
\\
\hline
\end{tabular}
\caption{The number in parenthesis indicates the multiplicity of the different solutions that occur
for the same number of Killing spinors.}
\end{center}
\end{table}

It is clear from table 2 that the backgrounds for which ${\rm hol}(\hat\nabla)\subseteq K\ltimes\bR^8$ can be characterized
by the number $L$ of parallel spinors, and the number $N$ of Killing spinors. This is because each case that appears has multiplicity
one.  This is not the case for backgrounds for which ${\rm hol}(\hat\nabla)$ is compact. For these some information about the embedding
of ${\cal K}_\tN$ in ${\cal P}_\tL$ is necessary to characterize the geometry.

Another result that becomes evident from table 2 is that, apart from the case with ${\rm Stab}(\e_1, \dots, \e_\tL)=\{1\}$,
for any given $L$, the Killing spinor equations have solutions
for any $1\leq N\leq L$. This is a consequence of the dilatino Killing spinor equation. However not all cases are independent.
For example, given ${\cal P}_\tL$, it is clear that all backgrounds with $N<L$ have the same $\hat\nabla$-parallel spinors. Therefore, one
expects that the geometry of all these backgrounds, called descendants in \cite{het2}, must have some similarity with
that of  backgrounds with $N=L$. This indeed is the case and it has been shown in \cite{het2} that a relation can be established
using  the  field equations of the theory. This relates the backgrounds lying horizontally in table 2.

There is also  a relation between the geometries
of backgrounds  lying   diagonally  in table 2. This will be described separately for the compact and non-compact cases   below.

\section{Spacetime geometry}

\subsection{Non-compact holonomy}

In table 2, there are  29 different types of  supersymmetric backgrounds for which
${\rm hol}(\hat\nabla)$ is {\it non-compact}. However it is not necessary to investigate them separately because
some of them are special cases of others.  This follows from the results of \cite{het2},
where all the Killing spinors are stated explicitly. To outline this relation consider the case of $(L,N)$ background, $N\not=7$,  i.e.~a background with $L$ parallel and $N$ Killing spinors, $N<L$.  As has already been mentioned the pair $(L,N)$ uniquely determines the background.
It turns out that the Killing spinors of this background are {\it identical}  to those of  $(N,N)$ background.
Thus the geometry of the $(L,N)$ backgrounds is a special case of that of $(N,N)$ backgrounds, $N\not=7$. This establishes
a relation between the geometries of backgrounds lying diagonally in table 2.
Therefore, it suffices to investigate the geometry of backgrounds for which all parallel spinors are Killing, i.e.~only that of the $(L,L)$
backgrounds for $L=1,2,3,4,5,6,8$. The $(8,7)$ backgrounds are special and should be treated separately. This has been done in \cite{het2}
and we shall not expand on this here.

\subsubsection{Geometry}

The spacetime  of $(L,L)$ backgrounds admits $\hat\nabla$-parallel or fundamental forms (\ref{forms}) constructed
from Killing spinor bilinears. It turns out that the fundamental forms of backgrounds with ${\rm hol}(\hat\nabla)=K\times \bR^8$ are
\bea
e^-~,~~~e^-\wedge \tau~,
\la{fforms}
\eea
where $e^-$ is a null 1-form and $\tau$ is a fundamental form of $K$.

The solution of the Killing spinor equations
\begin{itemize}

\item expresses  the 3-form $H$ in terms of the metric and fundamental  forms (\ref{fforms}), and

\item imposes restrictions on the geometry of spacetime which can be written as conditions
on the metric and (\ref{fforms}).

\end{itemize}

To describe both types of conditions  in detail, it is convenient to define the directions  ``transverse'' to the lightcone.
For this define the vector field $e_+$ using $ e^-(\cdot)=g(e_+, \cdot)$.
Since $e_+$ is also $\hat\nabla$-parallel, it spans a trivial bundle $I$ in $T\sp$. Moreover, one has
\bea
0\rightarrow I\rightarrow  {\rm Ker}\, \,e^-\rightarrow \xi_{T\sp}\rightarrow 0~,
\eea
where ${\rm Ker}\, e^-$ is spanned  by the vector fields $X$ of $\sp$ annihilated by $e^-$, $e^-(X)=0$. It is clear that
$\xi_{T\sp}$ has rank 8 and it is identified with the directions transverse to the lightcone.

In practise this means that one can adapt a local frame $(e^-, e^+, e^i)$, $i=1,\dots, 8$, where $e^+, e^i$ are defined
up to shifts along $e^-$, such that the solution of
the Killing spinor equations can be written as
\bea
ds^2&=&2 e^- e^++\delta_{ij}\,\, e^i e^j~,
\cr
H&=&e^+\wedge de^-+{1\over2}(h^{\mathfrak {k}}+h^{\mathfrak {k}^\perp})_{ij}\,\, e^-\wedge e^i \wedge e^j+ \tilde H~,
\la{solflux}
\eea
where
\bea
\tilde H={1\over3!}H_{ijk} e^i\wedge e^j\wedge e^k~.
\eea
We have already expressed some of the components of $H$ in terms of fundamental forms because these are universal.
 To identify the rest of the
components, first observe that the Lie algebra $\mathfrak{k}$ of $K$ is a subspace of $\Lambda^2(\bR^8)$, $\mathfrak{k}\subset \Lambda^2(\bR^8)$.
So one can use the metric to write $\Lambda^2(\bR^8)=\mathfrak{k}\oplus \mathfrak{k}^\perp$. So $h^{\mathfrak {k}}$ and $h^{\mathfrak {k}^\perp}$
are the components of the 2-form $h$ along $\mathfrak{k}$ and $\mathfrak{k}^\perp$, respectively. $\tilde H$ are the components of $H$
along the directions transverse to the lightcone. From now on forms denoted by tilde have components
only along the directions transverse to the light-cone.

The Killing spinor equations determine all components of $H$ apart from $h^{\mathfrak {k}}$. In particular, $h^{\mathfrak{k}^\perp}$
and $\tilde H$ are determined in terms of the metric and the fundamental forms. However these expressions
are case dependent. We shall mostly focus on  $\tilde H$. The expression for $h^{\mathfrak{k}^\perp}$ can be found in \cite{het1, het2}.

Next observe that
\bea
\hat\nabla e^-=0\Longleftrightarrow e_+~{\rm Killing}~,~~~de^-=e_+H~.
\eea
So $\sp$ admits a single null Killing vector field. This condition on the geometry is universal. There are additional conditions
which are case dependent. We shall mention these in the appropriate section.

\subsubsection{$Spin(7)\ltimes \bR^8$}
Let $\phi=\tilde \phi$ be the self-dual fundamental 4-form of $Spin(7)$.
In addition to the conditions that are universal and mentioned already, the Killing spinor equations  imply that
\bea
\tilde H=-\star \tilde d\phi+\star (\tilde\theta_\phi\wedge \phi)~,
\eea
and
\bea
&&\partial_+\Phi=0~,~~~
de^-\in \mathfrak{spin}(7)\oplus_s\bR^8~,~~~
\cr
&&2\partial_i\Phi-(\tilde\theta_\phi)_i-H_{-+i}=0~,~~~
\la{solspin7}
\eea
where
\bea
\tilde \theta_\phi=-{1\over6}\star(\star \tilde d\phi\wedge \phi)
\eea
is a Lee form, $\star$ is the Hodge duality operation along the transverse directions, and $\tilde d$ is the exterior
derivative again evaluated along the transverse directions. It is clear that $\tilde H$ can be expressed in terms of the fundamental
from $\phi$.  The expression is similar to that for 8-manifolds with a $Spin(7)$ structure and compatible $Spin(7)$
connection with skew-symmetric torsion \cite{stefan1}. The dilaton is invariant under the action of the vector field $e_+$.
The second condition in (\ref{solspin7}) is a geometric condition
which restricts the twist of the vector field $e_+$. In turn it implies that $e^-\wedge \phi$ is invariant under the action of $e_+$.
 The last condition can also be perceived
as a geometric condition which expresses the Lee form $\tilde \theta_\phi$ in terms of the dilaton.

\subsubsection{$SU(4)\ltimes \bR^8$}

Let $\omega_I=\tilde \omega_I$ and $\chi=\tilde \chi$ be the Hermitian  and the (4,0) fundamental forms of $SU(4)$, respectively.
$I$ is an almost complex structure in $\xi_{T\sp}$ associated with $\omega_I$. The  Killing spinor equations imply that
\bea
\tilde H=-i_{\tilde I}d\omega_I=-\star( \tilde d\omega_I\wedge \omega_I)-{1\over2}\star (\tilde\theta_{\omega_I}\wedge \omega_I\wedge \omega_I)~,
\la{hsu4}
\eea
and
\bea
&&\partial_+\Phi=0~,~~ de^-\in \mathfrak{su}(4)\oplus_s\bR^8~,~~~
\cr
&&\tilde {\cal N}(I)=0~,~~~~\tilde \theta_{\omega_I}=\tilde \theta_{{\rm Re}\,\chi}~,~~
\cr
&& 2\partial_i\Phi-(\tilde\theta_{\omega_I})_i-H_{-+i}=0~,
\la{solsu4}
\eea
where $\tilde{\cal N}$ is the Nijenhuis tensor of $I$ restricted along the
transverse directions and
\bea
\tilde \theta_{\omega_I}=-\star(\star \tilde d\omega_I\wedge \omega_I)~,~~~
\tilde\theta_{{\rm Re}\,\chi}=-{1\over4}\star(\star \tilde d{\rm Re}\,\chi\wedge {\rm Re}\,\chi)~,
\la{hlee}
\eea
are the Lee forms of $\omega_I$ and ${\rm Re}\,\chi$, respectively. The expression for $\tilde H$ is as that for the skew-symmetric torsion of the
 Bismut connection  for 2n-manifolds with a $U(n)$ structure, see also \cite{hull}-\cite{lust}.

There are two new type of conditions that appear in (\ref{solsu4}) compared to those which we have analyzed for the $Spin(7)\ltimes\bR^8$ case.
The first
is the vanishing of the Nijenhuis tensor for $I$. This is a consequence
of the dilatino Killing spinor equation.  The other is the equality between
the Lee forms $\tilde \theta_{\omega_I}$ and $\tilde \theta_{{\rm Re}\,\chi}$. This is required for the existence of a compatible
connection with skew-symmetric torsion on 8-dimensional manifolds with an $SU(4)$ structure.

\subsubsection{$Sp(2)\ltimes \bR^8$, $\times^2 SU(2)\ltimes \bR^8$, $SU(2)\ltimes \bR^8$ and $U(1)\times \bR^8$}

The fundamental forms of all these backgrounds are
\bea
e^-~,~~~~e^-\wedge \omega_r~,
\eea
where
$\omega_r=\tilde\omega_r$ are Hermitian forms on the space $\xi_{T\sp}$ transverse to the light-cone. These
can be thought of as the fundamental forms of the maximal compact subgroup $K$ in ${\rm Stab}(\e_1, \dots, \e_L)\equiv  K\times \bR^8$.
These Hermitian forms $\omega_r$ and their associated  endomorphisms $I_r$ have been explicitly given  in \cite{het2}.
The data provided can be re-organized more efficiently in terms of Clifford algebras. In particular using the
results of \cite{het2}, one can show that the typical fibre of $\xi_{T\sp}$ is an appropriate Clifford module as
indicated in table 3.

To see how the fundamental Hermitian forms can be identified from table 3, first consider the $SU(4)$ case.
It is clear that the almost complex structure $I$ can be thought of as the basis element of $\cl(\bR)$. Similarly,
it is known that the fundamental forms of $Sp(2)$ are Hermitian forms $\omega_r$  associated with an almost
hyper-complex structure $I_r$. Two of the almost complex structures, say $I_1$ and $I_2$, can be identified with the
two basis elements of $\cl(\bR^2)$. The third $I_3$ is the product of the other two, $I_3=I_1 I_2$, and so it is
represented by the even element of  $\cl(\bR^2)$ which again is the product of the two basis elements. This construction
is easily extended to all other cases. Note in particular that the geometry of both the $(8,7)$ and $(8,8)$ backgrounds can be described in this
way.

\begin{table}
\begin{center}
\begin{tabular}{|c|c|c|}
\hline
$N$ & ${\mathrm Stab}(\e_1,\dots,\e_L)$ & ${\rm Clifford}$
 \\
 \hline\hline
$2$ &  $SU(4)\ltimes\bR^8$&$\cl(\bR)$
\\ \hline
$3$ & $Sp(2)\ltimes\bR^8$&$\cl(\bR^2)$
\\ \hline
$4$ & $(\times^2 SU(2))\ltimes\bR^8$&$\cl(\bR^3)$
\\ \hline
$5$ & $SU(2)\ltimes\bR^8$&$\cl(\bR^4)$
\\ \hline
$6$ & $U(1)\ltimes\bR^8$&$\cl(\bR^5)$
\\ \hline
$7$ & $\bR^8$&$\cl(\bR^6)$
\\ \hline
$8$ & $\bR^8$&$\cl(\bR^7)$
\\
\hline
\end{tabular}
\caption{The number of Killing spinors is given in the first column. In the second column the isotropy group
of the parallel spinors is given. In the last column the associated Clifford algebra of endomorphisms is indicated.}
\end{center}
\end{table}

Now $\tilde H$ can be given as in (\ref{hsu4}) with respect to any of the endomorphisms $I_r$, say $I=I_1$.
The rest of the conditions  implied from the  Killing spinor equations are
\bea
&&\partial_+\Phi=0~,~~~i_{ I_r}d\omega_r=i_{ I_s}d\omega_s~,~~~r\not=s
\cr
&&de^-\in \mathfrak{k}\oplus_s\bR^8~,~~~ \tilde {\cal N}(I_r)=0~,~~~~
\tilde \theta_r=\tilde \theta_s~,~~~r\not=s
\cr
&&2\partial_i\Phi-(\tilde\theta_r)_i-H_{-+i}=0~,~
\eea
where $\tilde\theta_r$ is the Lee form of $\omega_r$ given as in (\ref{hlee}). The above  conditions
can be easily derived from those of the Killing spinor equations for the  $SU(4)\ltimes\bR^8$ case. This can be done
by requiring that the conditions that are valid for the $I$ endomorphism should now be valid for
all $I_r$ endomorphisms.

\subsubsection{$\bR^8$} \la{r8}

It remains to state the conditions for the (8,8) case. It turns out that the Killing spinor equations imply that
\bea
e^-\wedge de^-=0~,~~~\tilde H=0~.
\la{conr8}
\eea
In section \ref{half}, we shall classify all such backgrounds by solving both the above conditions and associated field equations.

\subsection{Compact holonomy}

Table 2 indicates that there are 32 cases that we should consider. So it is natural to seek a simplification similar to that
we have  introduced for the non-compact cases. It turns out that there is a simplification but not as effective to reduce
the analysis as for the non-compact holonomy cases. This is because the geometry depends on the embedding of ${\cal K}_N$ in ${\cal P}_L$.
Therefore more information is needed to determine the geometry than  just the
dimension of these spaces.

To give an example where a simplification can be made, consider the two distinct cases that arise in  $N=2$ backgrounds with
${\rm hol}(\hat\nabla)\subseteq SU(3)$. Inspecting the Killing spinors in \cite{het2}, it is easy to see that one of these two cases is a special case
of $N=2$ backgrounds with ${\rm hol}(\hat\nabla)\subseteq SU(4)\ltimes \bR^8$,  and the other case is a special case
of $N=2$ backgrounds with ${\rm hol}(\hat\nabla)\subseteq G_2$. However there are several cases that occur in backgrounds with
${\rm hol}(\hat\nabla)\subseteq SU(2)$ which do not have such an association. Because of this, we shall describe the geometry of
backgrounds for which all parallel spinors are Killing, i.e.~$N=L$. The case with $L=16$ corresponds to the maximally supersymmetric
backgrounds and it is known that these are locally isometric to $\bR^{9,1}$ \cite{gpjose}.

\subsubsection{Geometry}\la{comgen}

The $\hat\nabla$-parallel forms which arise as Killing spinor bilinears are
\bea
e^a~,~~~\tau~,
\eea
where $e^a$ are 1-forms and $\tau$ are the fundamental forms of $K\equiv {\rm Stab}(\e_1,\dots\e_\tL)$, ${\rm hol}(\hat\nabla)\subseteq K$.
The number of parallel 1-forms
depends on $K$ and one of them is always time-like. The minimal number of parallel 1-forms are 3,4 and 6 for $G_2$, $SU(3)$ and $SU(2)$, respectively.

Let $e_a$ denote the dual vector field of $e^a$, $e^a(\cdot)=g(e_a, \cdot)$. Provided that $dH=0$, the commutator $[e_a, e_b]$ is again
a $\hat\nabla$-parallel vector field. The Killing spinor equations in most cases do not put sufficient
restrictions on the commutator $[e_a, e_b]$ to express it in terms of the original vector field $e_a$. So potentially $\sp$ may admit
more $\hat\nabla$-parallel vector fields than those constructed from Killing spinor bilinears.
To simplify the analysis that follows, we shall  assume that the vector fields constructed form Killing spinor
bilinears span a Lie algebra under Lie brackets. Assuming that the action of the vector fields can be integrated to
a free group action, $\sp$ is a principal bundle,
$\sp=P(G,B;\pi)$, where the fibre $G$ has Lie algebra that of the vector field $\{e_a\}$ and $B$ is the base space. Moreover $\sp$
is equipped with a principal bundle
connection $\lambda^a\equiv e^a$. The Lie algebras of the fibre groups $G$ have been tabulated in table 4.

\begin{table}
\begin{center}
\begin{tabular}{|c|c|c|}
\hline
${\mathrm Stab}(\e_1,\dots,\e_L)$ & $1-{\rm forms}$&$\mathfrak{Lie}\,G$
 \\
 \hline
$G_2$&$3$&$\bR^3~,~\mathfrak{sl}(2,\bR)$
\\
\hline
$SU(3)$ & $4$&$\bR^4~,~\mathfrak{sl}(2,\bR)\oplus \bR~,~ \mathfrak{su}(2)\oplus \bR~,~\mathfrak{cw}_4$
\\
\hline
$SU(2)$ & $6$&$ \bR^6~,~\mathfrak{sl}(2,\bR)\oplus\mathfrak{su}(2)~,~\mathfrak{cw}_6$
\\
\hline
\end{tabular}
\end{center}
\caption{In the first column, the compact isotropy groups of spinors are stated. In the second column, the number of
1-form spinor bilinear is given. In the third column, the associated Lorentzian Lie algebras are exhibited.
The structure constants of the 6-dimensional Lorentzian Lie algebras of the $SU(2)$ case
are self-dual.}
\end{table}

Using the principal bundle data, the solution of the  Killing spinor equations can be expressed as
\bea
ds^2&=&\eta_{ab}\, \lambda^a \lambda^b+\pi^* \tilde g
\cr
H&=&{1\over3} \eta_{ab} \lambda^a\wedge  d\lambda^b+{2\over3} \eta_{ab} \lambda^a\wedge {\cal F}^b+ \pi^* \tilde H~,
\la{gh}
\eea
where
\bea
{\cal F}^a\equiv d\lambda^a-{1\over2} H^a{}_{bc} \lambda^b\wedge \lambda^c~,
\eea
is the curvature of $\lambda$ and $\tilde g=\delta_{ij}\, e^i e^j$ is a metric on $B$. Apart from the above conditions which are universal,
the Killing spinor equations impose  additional restrictions on the geometry of spacetime which depend on $K$.  These will be given when
we describe each case separately. Observe that $H$ is the sum of the Chern-Simons form of $\lambda$ and a 3-form $\tilde H$ of $B$.
From now one, the forms on $B$ will be denoted with a tilde.

\subsubsection{$G_2$}

Let $e^a$, $a=0,1,2$ and $\varphi=\tilde\varphi$ be three 1-forms and the fundamental $G_2$ form, respectively. In addition to (\ref{gh}), the
Killing spinor equations imply
\bea
\tilde H=-{1\over6} (\tilde d\varphi, \star \varphi)\, \varphi+ \star \tilde d\varphi-\star (\tilde\theta\wedge\varphi)
\la{g2h}
\eea
and
\bea
&&\partial_a\Phi=0~,~~~{\cal F}\in \mathfrak{g}_2~,~~~\epsilon^{abc} H_{abc}+H_{ijk}\varphi^{ijk}=0~,
\cr
&&\tilde\theta=2d\Phi~,~~~\tilde d\star\varphi=-\tilde\theta\wedge \star\varphi~,
\la{g2con}
\eea
where
\bea
\tilde\theta=-{1\over3} \star(\star \tilde d\varphi\wedge \varphi)~,
\eea
is the Lee form of $\varphi$, and $\tilde d$ and  $\star$ is the exterior derivative  and  the Hodge operation on $B$, respectively.
It is clear that the dilaton is invariant under all $e_a$ vector field and so is a function of $B$. Moreover ${\cal F}\in \mathfrak{g}_2$
implies that the principal bundle connection is a $G_2$ instanton of $B$. Another consequence of the same condition is that $\varphi$ is
invariant and since $i_a\varphi=0$, it is the pull-back of a form on $B$. In fact the 7-dimensional manifold  $B$ admits a $G_2$ structure.
All the conditions in (\ref{g2con}) arise from the dilatino Killing spinor equation apart from the last one.  This is required \cite{stefan2}
for $(B, \tilde g, \tilde H)$ to admit a compatible  metric connection with skew-symmetric torsion, $\hat{\tilde \nabla}$, and holonomy contained in $G_2$.

The expression for $\tilde H$ depends on whether $G$ is abelian or not. If $G$ is abelian, then the first term in (\ref{g2h}) for
 $\tilde H$ vanishes  as can be seen from (\ref{g2con}). On the other hand if $G=SL(2,\bR)$, then the same term becomes proportional
 to the volume of $SL(2,\bR)$.

\subsubsection{$SU(3)$}
Let $e^a$, $a=0,1,2,3$, and $\omega=\tilde\omega$ and $\chi=\tilde\chi$ be four 1-forms, and the Hermitian and (3,0) fundamental forms of $SU(3)$, respectively.
In addition to (\ref{gh}),
 the Killing spinor equations imply that
\bea
\tilde H=-i_{\tilde I}\tilde d \omega=\star\tilde d\omega-\star(\tilde\theta_\omega\wedge \omega)~,
\eea
and
\bea
&&\partial_a\Phi=0~,~~~{1\over3!} \e^{abcd} H_{bcd}-{1\over2}{\cal F}^a_{ij}\omega^{ij}=0~,~~~~({\cal F}^a)^{2,0}=0~,
\cr
&&\tilde{\cal N}(I)=0~,~~~\tilde\theta_\omega=\tilde\theta_{{\rm Re}\chi}~,
\cr
&&\partial_i\Phi-{1\over2} (\tilde\theta_\omega)_i=0~,
\la{consu3}
\eea
where
\bea
\tilde\theta_\omega=-\star(\star\tilde d\tilde\omega\wedge \tilde\omega)~,~~~
\tilde\theta_{{\rm Re}\chi}=-{1\over2} \star(\star\tilde d{\rm Re}\chi\wedge {\rm Re}\chi)~,
\eea
are the Lee forms of $\omega$ and $\chi$, respectively. The conditions have similarities with those of both the $G_2$ and $SU(4)\ltimes\bR^8$
cases. The dilaton $\Phi$ is a function of $B$.

The geometry of $B$ depends on whether $G$ is abelian or non-abelian. If $G$ is abelian, then $(B,\tilde g, \tilde H)$ is a complex manifold that
admits a compatible metric connection, $\hat{\tilde \nabla}$,  with skew-symmetric torsion and with holonomy contained in $SU(3)$.
This follows from
\bea
 {\cal F}\in  \mathfrak{su}(3)~,
 \eea
 which in turn implies that $\omega$ and $\chi$ are invariant under the action of all vector field $e_a$,
  and the equality of Lee forms in (\ref{consu3}).
Moreover $\lambda$ is a Donaldson type of connection. Since $2d\Phi=\tilde\theta_\omega$,
$B$ is conformally balanced\footnote{ It is known that the smooth compact
$2n$-dimensional conformally balanced manifolds $B$ with ${\rm hol}(\hat{\tilde \nabla})\subseteq SU(n)$ and with  $\tilde d \tilde H=0$
 \cite{ivanovgp} are Calabi-Yau  with $\tilde H=0$. However, there
are non-compact smooth examples.}.

Next suppose that $G$ is non-abelian and so is locally either $\bR\times SU(2)$ or $SL(2,\bR)\times U(1)$.
In such a case,
\bea
{\cal  F}\in \mathfrak{su}(3)\oplus \bR~,
\eea
and $\chi$ is not invariant under the $\bR$ and $U(1)$ group actions, respectively. Instead it is invariant up to a $U(1)$ rotation. As a result, the
canonical bundle of $B$ is twisted and so $B$ has not an $SU(3)$ structure but rather a $U(3)$ one.
So in this case, $(B, \tilde g, \tilde H)$ is a Hermitian manifold with a compatible
connection, $\hat{\tilde \nabla}$,  with skew-symmetric torsion and with holonomy contained in $U(3)$.

\subsubsection{$SU(2)$} \la{su2}

Let $e^a$, $a=0,\dots,5$, and $\omega_r=\tilde\omega_r$, $r=1,2,3$,  be six 1-forms and the three Hermitian fundamental forms of $SU(2)$, respectively,
where the endomorphism $I_r I_s=-\delta_{rs} {\bf 1}_{4\times 4}+\epsilon_{rst} I_t$.
In addition to the conditions (\ref{gh}), the Killing spinor equations imply that
\bea
\tilde H=-i_{I_1}\tilde d\omega_1
\eea
and
\bea
&&\partial_a\Phi=0~,~~~~,~~~H_{a_1a_2a_3}+{1\over3!} \epsilon_{a_1a_2a_3}{}^{b_1b_2b_3} H_{b_1b_2b_3}=0~,
\cr
&&
i_{I_r}\tilde d\omega_s=i_{I_s}\tilde d\omega_r~,~~r\not=s~,~~~{\cal F}^a\in \mathfrak{su}(2)~,~~~2\partial_i \Phi-(\tilde\theta_{\omega_1})_i =0~,
\la{consu2}
\eea
where $\tilde\theta_{\omega_1}$ is the Lee form of $\omega_1$ as in (\ref{hlee}). Again $\Phi$ is a function of the base space $B$.
In this case, the Lie algebra of the fibre $G$ is self-dual, i.e.~the structure constants satisfy the self-duality condition. In addition,
 the principal bundle connection $\lambda$ is an anti-self-dual
instanton. It turns out that the conditions (\ref{consu2}) imply that the base space $B$ is {\it conformally hyper-K\"ahler}.

\section{All half supersymmetric solutions}

These are the solutions of both the Killing spinor and field equations that admit 8 Killing spinors. There are three
classes of such backgrounds which have been classified in \cite{het4}. One class is that of $N=L=8$ backgrounds with ${\rm Stab}(\e_1, \dots, \e_\tL)=SU(2)$ investigated in \ref{su2}.
The other class is that of $N=L=8$ backgrounds with ${\rm Stab}(\e_1, \dots, \e_\tL)=\bR^8$ examined in \ref{r8}. The third
class  is that of $(L,N)=(16, 8)$ backgrounds associated
with  ${\rm Stab}(\e_1, \dots, \e_{16})=\{1\}$. It turns out that the third case is a special case of the other two. So we have
only two possibilities to investigate.

\subsection{$SU(2)$}

It has been demonstrated that the spacetime in this case is a principal bundle $\sp=P(B,G;\pi)$ over a conformally hyper-K\"ahler manifold $B$,
 equipped with a anti-self-dual connection $\lambda$
and fibre group $G$ with a self-dual Lorentzian Lie algebra. Using these data, one can write
\bea
ds^2&=&\eta_{ab} \lambda^a \lambda^b+h \, ds_{\rm hk}^2~,~~~e^{2\Phi}=h~,
\cr
H&=&{1\over3} \eta_{ab} \lambda^a\wedge d\lambda^b+{2\over3} \eta_{ab} \lambda^a\wedge {\cal F}^b-\star_{\rm hk} \tilde dh~,
\eea
where $h$ is a function of $B$ and $ds_{\rm hk}^2$ is a 4-dimensional hyper-K\"ahler metric.

To find explicit examples one has to specify a 4-dimensional hyper-K\"ahler manifold, a anti-self dual instanton connection over it
with gauge group $G$
and to determine the function $h$. The latter  is found by exploring the Bianchi identity (\ref{anb}) of $H$, i.e.~$dH=0$, where
we have neglected the anomaly term which is proportional to $\a'$. In particular,
\bea
dH\equiv d\pi^*\tilde H+\eta_{ab} {\cal F}^a\wedge {\cal F}^b=0~.
\la{closeH}
\eea
One can include higher order corrections $\a'$ corrections and the complete analysis has been done \cite{het4}.
This in turn can be written as
\bea
-\nabla^2_{\rm hk} h-{1\over2} \eta_{ab}\,{\cal F}^a_{ij}\,\, {\cal F}^{bij}=0~,
\la{fineqn}
\eea

One class of solutions is given by taking $\lambda$ to be a trivial connection. In such a case the spacetime is
${\sp}=G\times B$ and
\bea
ds^2=ds^2(G)+ds^2(B)~,~~~H={1\over6} H_{abc} \lambda^a\wedge \lambda^b\wedge \lambda^c~,~~~e^{2\Phi}={\rm const}~,
\eea
i.e.~$H$ is determined in terms of the structure constants of $G$.

An example of a solution with non-trivial connection $\lambda$ can be constructed by taking $B=\bR^4$ and $G=SL(2,\bR)\times SU(2)$.
In addition consider an anti-self dual connection $\lambda$ on $\bR^4$ with gauge group $SU(2)$ and  with instanton number 1.
Since only the $SU(2)$ subgroup of $G$ is gauged, the spacetime is $\sp=SL(2,\bR)\times X_7$. In particular,
\bea
ds^2&=&ds^2\big(SL(2, \bR)\big)+\delta_{pq} \lambda^p \lambda^q+h \, ds(\bR^4)~,~~~e^{2\Phi}=h~,~~~p,q=1,2,3~,
\cr
H&=&d{\rm vol}(SL(2,\bR)+{1\over3} \delta_{pq} \lambda^p\wedge d\lambda^q+{2\over3} \delta_{pq} \lambda^p\wedge {\cal F}^q-\star_{\rm hk} \tilde dh~,
\eea
where
\bea
h=1+4 {|x|^2+2\rho^2\over (|x|^2+\rho^2)^2}~,~~~~~x\in \bR^4~,
\eea
and where $\rho$ is the size of the instanton. This solution easily generalizes to multi-instanton $SU(2)$ solutions \cite{het4}.
Thus there is a class of solutions which depends on $8\nu-3$ parameters, the moduli of $SU(2)$ instantons with instanton number $\nu$.

\subsection{$\bR^8$} \la{half}

The conditions stated in (\ref{conr8}) for this case imply that there is a choice of coordinates $(u,v, x^i)$ such that

\bea
ds^2=2 e^- e^++ ds^2(\bR^8)~,~~~~H= d(e^-\wedge e^+)~,~~~
\cr
e^-=h^{-1} dv~,~~~e^+=du+V dv+ n_i dx^i~.
\eea
All components of the metric and $H$ depend on $v$ and $x$, and $e_+=\partial_u$ is the null parallel vector field.

The solutions of the Killing spinor equations are determined up to the functions $h$ and $V$, and the 1-form $n$.
These in turn can be found by solving the field equations (\ref{feqn}). If in addition one assumes that
$h$, $V$ and $n$ are $v$ independent, then the field equations imply that
\bea
\partial_i^2 h=\partial_i^2 V=0~,~~~\partial^i dn_{ij}=0~.
\eea
So $h$ and $V$ are harmonic functions of $\bR^8$ and $dn$ satisfies the Maxwell equations on $\bR^8$.
The solution is a superposition of fundamental strings \cite{rr}, pp-waves and  null rotations.

\vskip 0.5cm
\noindent{\bf\large Acknowledgements} \vskip 0.3cm
We would like to thank J.~Gutowski and D.~Roest for many useful discussions.

\end{document}